\documentclass[leqno,draft]{article}

\begin{document}

\title{Some remarks about Cauchy integrals \\ and fractal sets}

\author{Stephen Semmes \\
        Rice University}

\date{}

\maketitle

        If $\mu$ is a finite Borel measure on the complex plane, then
the Cauchy integral
\begin{equation}
        C(\mu)(z) = \int_{\bf C} \frac{1}{z - \zeta} \, d\mu(\zeta)
\end{equation}
defines a holomorphic function of $z$ on the complement of the support
of $\mu$.  For simplicity, let us restrict our attention to measures
with compact support in ${\bf C}$, although one can also make sense of
Cauchy integrals of measures with noncompact support and infinite mass
under suitable conditions.  In the classical situation where $\mu$ is
supported on a nice curve, the Cauchy integral of $\mu$ has a jump
discontinuity across the curve.  For more regular measures, the Cauchy
integral converges absolutely for every $z \in {\bf C}$ and defines a
continuous function on the plane.  Of course, this function is not
holomorphic in any neighborhood of the support of $\mu$.  For any
finite measure $\mu$ on ${\bf C}$, the Cauchy integral $C(\mu)(z)$
makes sense as a locally integrable function on ${\bf C}$, whose
$\overline{\partial}$ derivative is a constant multiple of $\mu$ in
the sense of distributions.

        On nice regions in ${\bf C}$, the Cauchy integral formula can
be used to recover arbitrary holomorphic functions from their boundary
values, under suitable conditions.  The Cauchy integral leads to a
projection from general functions on the boundary to boundary values
of holomorphic functions.

        In more fractal situations, one might prefer to think of
Cauchy integrals as a way of solving $\overline{\partial}$ problems,
to make corrections to get holomorphic functions instead of producing
them directly.  One might view this as being more like several complex
variables.  A basic scenario would be to multiply a holomorphic
function by a non-holomorphic function with some regularity, and to
try to make some relatively small corrections to the product to get a
holomorphic function.  Perhaps the holomorphic function has nice
boundary values and the other function is defined on the boundary, and
the correction is intended to yield the boundary values of a
holomorphic function.  For example, the product of a holomorphic
function and a rational function with poles in the interior may not be
holomorphic, and this can be corrected with a finite-rank operator to
get rid of the singularities.

        Bergman spaces and projections could also be considered on the
interior.  It would still be regularity of the non-holomorphic
functions up to the boundary that matters, though.

        In ${\bf R}^n$, one can use Cauchy integrals associated to
Clifford analysis.  However, remember that the product of Clifford
holomorphic functions is not necessarily holomorphic, because
of noncommutativity of the Clifford algebra,

        In any dimension, it is already somewhat interesting to look
at Bergman spaces of locally constant functions on regions with
infinitely many connected components.

        There seem to be a lot of choices involved with fractals.
This may be a bit disconcerting, but perhaps it is also just as well.

\end{document}